\documentclass[12pt]{amsart}
\usepackage{floatflt,epsfig,subfigure}
\usepackage{a4}
\usepackage{version}

\usepackage{amsmath}
\usepackage{amssymb}

\def\proofof #1 {\par\medskip\noindent {\sc Proof of #1. }}

\def\sketchof #1 {\par\medskip\noindent {\sc Sketch of proof of #1. }}

\newcommand{\N}{\mathbb{N}}

\newcommand{\Z}{\mathbb{Z}}

\newcommand{\C}{\ensuremath{\mathbb{C}}}

\newcommand{\D}{\mathbb{D}}

\renewcommand{\H}{\mathbb{H}}

\newcommand{\re}{\operatorname{Re}}
\newcommand{\im}{\operatorname{Im}}

\newcommand{\cl}[1]{\overline{#1}}

\newcommand{\Ek}{E_{\kappa}}

\newtheorem*{thm}{Theorem}
\newtheorem{prop}{Proposition}

\newcommand{\PC}{\mathcal{P}}

\def\authoraddress#1{
    \def\thefootnote{\fnsymbol{footnote}}
	  \footnotetext[0]{\hskip -4ex{\bf Address:} #1\\}
    \def\thefootnote{\arabic{footnote}}
}

\title[An Answer to a Question of Herman, Baker and Rippon]{An Answer to a Question of Herman, Baker and Rippon 
   Concerning Siegel Disks}

\author{Lasse Rempe}

\begin{document}

 \begin{abstract}
  Consider the family of exponential maps $\Ek(z)=\exp(z)+\kappa$. 
  This article shows that any unbounded Siegel disk $U$ of $\Ek$ contains
  the singular value $\kappa$ on its boundary.
  By a result of Herman, this
  implies that $\kappa\in \partial U$ if the rotation
  number is diophantine.
 \end{abstract}

 \maketitle

\authoraddress{Lasse Rempe, Mathematisches Seminar der CAU Kiel,
  Ludewig-Meyn-Str. 4, 24098 Kiel, Germany, 
  {\tt lasse$@$math.uni-kiel.de}}

 \section{Introduction}
 In the collection \cite{article:haymanlist} of research problems in complex
 analysis, the following was posed as problem
 2.86, attributed to Herman, Baker and Rippon. 
 Let the function $f_{\lambda}(z)=\lambda(e^z-1)$, $|\lambda|=1$, have a
 Siegel disk $U$ that contains $0$.
  \begin{enumerate}
    \item[(a)] Prove that there exists some number $\lambda$ such that
               $U$ is bounded in $\C$.
    \item[(b)] If $U$ is unbounded in $\C$, does the singular value
      $-\lambda$ belong to $\partial U$?
  \end{enumerate}
 In \cite{article:herman}, it was shown that if $\lambda=e^{2\pi i
 \alpha}$, where $\alpha$ is diophantine, then $U$ is unbounded.
 Rippon \cite{article:rippon} generalized an argument of Carleson and
 Jones \cite[page~86]{book:carlesongamelin}
 to give an elementary
 proof that $-\lambda\in\partial U$ for almost every $\lambda$.
 It was also mentioned in \cite{article:rippon} that
 (a) could be solved by adapting a method from
 \cite{article:douady}.

 In this note, we give a positive answer to (b), with a rather simple
 proof. 
  This implies in
 particular that the singular value lies on the boundary for
 diophantine rotation numbers. We will prove the result in the
 following form, which allows Siegel disks of arbitrary period.
 \begin{thm}
  \label{thm:main}
  Let $\kappa\in\C$ and suppose that the function 
   $E_{\kappa}(z):= \exp(z)+\kappa$ 
   has an unbounded Siegel disk $U$. Then there
   is $j$ such that $\kappa\in \partial \Ek^j(U)$.
 \end{thm}
 Note that $f_{\lambda}$ is conjugate to $E_{\kappa}$ for
   $\kappa=\log\lambda - \lambda$.

\hspace{0.5cm}
\textsc{Acknowledgements} This work was conducted during a visit
 to the Institute for
 Mathematical Sciences at the State University of New York at Stony
 Brook, for whose continuing support and hospitality I am most grateful.
 I also thank Walter Bergweiler, Rodrigo
 Perez and Saeed Zakeri for interesting discussions on
 this work.

\section{Proof of the Theorem}

  For basic definitions and results,
   we refer the reader to an expository text on the iteration of
   entire functions such as \cite{article:baker,article:walter}.

 To prove the theorem, we shall use the following
 property of exponential
 maps:
 \begin{prop} \label{prop:curveinjuliaset}
  Let $\kappa\in\C$. Then 
   there exists a curve $\gamma:[0,\infty)\to J(\Ek)$ such that
  $\lim_{t\to\infty} \re(\gamma(t)) = +\infty$ and
  $\limsup_{t\to\infty} |\im(\gamma(t))| < \infty$.
 \end{prop}
 This well-known fact seems to have been
  first proved by Devaney, Goldberg and Hubbard
  \cite{preprint:dgh}. In fact, it is
  now known that the set of escaping points of $\Ek$ consists entirely
  of such curves \cite{article:dierkesc}.

  Fix a $\kappa\in\C$ for which $\Ek$ has a Siegel disk
 $U$. Sullivan's result that rational functions do not have wandering
 domains \cite{article:sullivan}
 has been generalized to the family of exponential maps by
 Baker and Rippon \cite{article:bakerrippon}. It is well known that
  \[ \partial U \subset \PC := \cl{\{\Ek^n(\kappa):n\in\N\}}; 
       \tag{*} \label{eqn:bdy} \]
 (see
 e.g. \cite[Theorem 7]{article:walter}).
 Thus $\kappa$ belongs to the Julia set of $\Ek$. Indeed,
 otherwise $\kappa$ would eventually map into some periodic component
 of the Fatou set. Since each such component is either an attracting
 domain, a parabolic domain, a Baker domain\footnote{%
 In fact, Baker domains do not occur for exponential maps as all
 escaping points lie in the Julia set
 \cite{article:alexmisha}.}
  or a Siegel disk
 \cite[Theorem 6]{article:walter}, $\PC$
 would then intersect the Julia set in at most finitely many
 points, which contradicts (\ref{eqn:bdy}).

 Let us suppose that $\kappa\notin \partial
 \Ek^j(U)$ for every $j\in\N$; we wish to show that $U$ is bounded. 
 Choose $\delta>0$ such that 
  \[\cl{\D_{\delta}}(\kappa) :=
      \{z\in\C:|z-\kappa|\leq \delta\} 
       \subset \C\setminus\bigcup_j \Ek^j(U). \]
 Note that this implies that, for every $j$,
  $\Ek^j(U) \subset \C \setminus \cl{\H}$,
 where 
  $\H=\Ek^{-1}\bigl(\D_{\delta}(\kappa)\bigr)
          =\{z\in\C:\re z < \log \delta\}$.

 \begin{prop}
  For every $j\in\N$, the set $\Ek^j(U)$ has bounded imaginary part.
 \end{prop}
 \begin{proof}
   Let $\gamma$ be the curve given by Proposition
  \ref{prop:curveinjuliaset}.
  Since $\kappa\in J(\Ek)$, we can find a preimage $g_0$ of $\gamma(0)$
 under an iterate of $\Ek$ such that $|g_0 - \kappa| < \delta$. 
 Taking the appropriate pullback of $\gamma$, we obtain a curve
 $g:[0,\infty)\to J(\Ek)$ with $g(0)=g_0$. (If the orbit of $\kappa$
 intersects $\gamma$, we might not be able to take these
 pullbacks. However, in this case we can pull back $\gamma$ along the
 orbit of $\kappa$ and obtain a curve actually starting at $\kappa$.)
 Choose the largest $t_0$ with
 $|g(t_0)-\kappa| = \delta$ and consider the set
  \[ K := \Ek^{-1}\left(\cl{\D_{\delta}(\kappa)} 
                          \cup g(\,[t_0,\infty)\,)\ \right). \]

\begin{figure}
 \center
  \epsfig{file=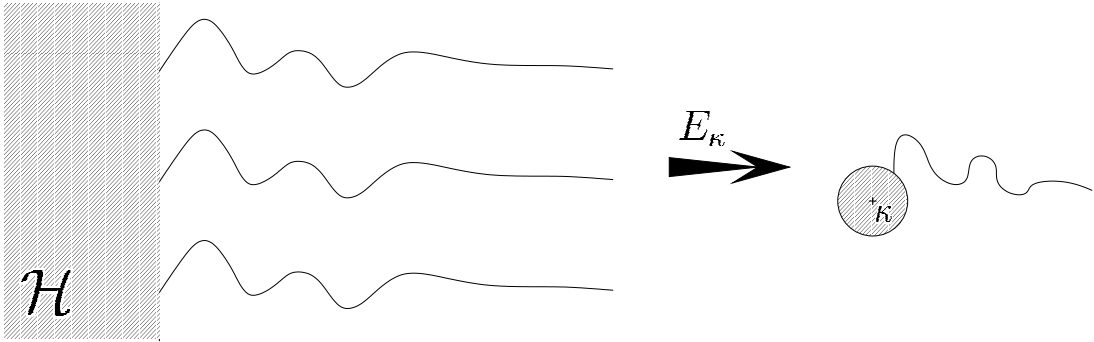, width=\textwidth}
 \caption{ \label{figK} The set $K$ and its image.}
\end{figure}%

 This set consists of $\cl{\H}$ together with the preimages of
 $g\bigl([t_0,\infty)\bigr)$ (see Figure~\ref{figK}). 
 Each of these preimages is asymptotic to a
 line $\{\im z = 2k\pi\}$, where $k\in\Z$, 
 and thus has bounded imaginary part. Therefore
 every component of $\C\setminus K$ has bounded imaginary part. Since
 $\Ek^j(U) \subset \C\setminus K$, this proves the claim.
\end{proof}

\begin{proof}[of the Theorem]
  By the previous proposition, we can find $S>0$ with
  $|\im z|<S$ for all $j\in\N$ and $z\in \Ek^j(U)$. Choose $R>1$
  large enough such that $\exp(R-1) > -\log\delta$ and
   $\exp(R)\geq \exp(R-1) + S + 1 + |\kappa|$.
  Then, if $z\in\C$ with $r := \re z \geq R$
   and $|\im\Ek(z)| < S$,
   \[
     |\re\Ek(z)| > |\Ek(z)| - S \geq
         \exp(r) - |\kappa| - S
         \geq  \exp(r-1) + 1. 
   \]
  If also $\re \Ek(z) > \log\delta$,
    then in particular $\re\bigl(\Ek(z)\bigr) - 1 > \exp(r-1)$. 
  
  Now suppose that there was a $z\in U$ with
   $r=\re(z)\geq R$. It then follows by
   induction that
   \[ \re\bigl(\Ek^n(z)\bigr) - 1 > \exp^n(r-1), \]
   which is a contradiction since points in a Siegel disk do not
   escape to $\infty$.
 \end{proof}

\end{document}